\documentclass[3p,authoryear,12pt]{elsarticlespe}
\usepackage{a4wide}

\usepackage{amsmath,amssymb,amsthm}

\newtheorem{thm}{Theorem}
\newtheorem{coro}{Corollary}
\newtheorem{assh}{Hypothesis}
\newdefinition{defi}{Definition}
\newenvironment{mypr}
{
\par\vspace{.5\baselineskip}\noindent\textbf{Proof.}}
{\nobreak\hfill$\Box$\par\vspace{.5\baselineskip}}

\usepackage{ifthen}

\usepackage{algorithm,algpseudocode}

\usepackage{color,pstcol,pstricks}
\newgray{backgroundgray}{0.90}
\definecolor{backgroundgray}{gray}{0.90}

\newenvironment{petit}
{\par\vspace{.5\baselineskip}\noindent\footnotesize}
{\nobreak\par\vspace{.5\baselineskip}}

\def\tempfolder{./temp/}
\providecommand{\imprimeTexteCache}{non}
\providecommand{\imprimeTexteSecret}{non}

\newif\ifsol
\newif\ifnfs

\ifthenelse{\equal{\imprimeTexteCache}{oui}}{\soltrue}{\solfalse }
\ifthenelse{\equal{\imprimeTexteSecret}{oui}}{\nfstrue}{\nfsfalse}

\ifsol
   \providecommand{\sol}[1]{\bs #1 \es}
\else
   \providecommand{\sol}[1]{}
\fi

\ifnfs
   \providecommand{\nfs}[1]{~\\{\footnotesize \sffamily \emph{\textbf{Note: } #1}~\\}}
\else
   \providecommand{\nfs}[1]{}
\fi

\newenvironment{mySolution}
{\red \it
\par\vspace{.5\baselineskip}\noindent\textbf{Solution.~}}
{\vspace{.5\baselineskip}}
\def\bs{\begin{mySolution}}
\def\es{\end{mySolution}}

\def\problemfile#1{%
\begin{filecontents}{\tempfolder#1}}

\ifnfs

\else

\fi

\ifnfs

\else

\fi

\ifnfs

\else

\fi

\ifnfs
\newcommand{\tbd}[1]{\textsc{To be done:\\#1}}
\else
\newcommand{\tbd}[1]{}
\fi

\newtheorem{maquestion}{\sc Question}[section]

\ifsol

\else

\fi

\newcommand{\imp}[1]{{\blue \bf #1}}

\newcommand{\bp}{ % begin program
  \small \ttfamily
 \begin{tabbing}
 aaa\=aaa\=aaa\=aaa\=aaaaaaaa \= aaaaaaaaaa\= \kill
 }
\newcommand{\ep}{\end{tabbing}\normalfont\normalsize }% end program

\newcommand{\sref}[1]{Section~\ref{#1}}
\newcommand{\cref}[1]{Chapter~\ref{#1}}
\newcommand{\coref}[1]{Corollary~\ref{#1}}
\newcommand{\eref}[1]{Eq.(\ref{#1})}

\newcommand{\thref}[1]{Theorem~\ref{#1}}

\newcommand{\calC}{\mathcal C}

\newcommand{\calF}{\mathcal F}

\newcommand{\calP}{\mathcal P}

\newcommand{\lp}{\left(}
\newcommand{\lb}{\left[}
\newcommand{\lc}{\left\{}

\newcommand{\rp}{\right)}
\newcommand{\rb}{\right]}
\newcommand{\rc}{\right\}}

\newcommand{\eps}{\epsilon}
\newcommand{\esp}[1]{\E\left(#1\right)}
\newcommand{\espc}[2]{\E\left(\left.#1 \right| #2\right)}

\newcommand{\abs}[1]{\left| #1 \right|}

\newcommand{\nin}{\not\in} %% negation of \in
\newcommand{\ind}[1]{\mathbf{1}_{\{#1\}}}

 \newcommand{\norm}[1]{\left\|#1\right\|}

\def\be{\begin{equation}}
\def\ee{\end{equation}}
\def\ben{\[}
\def\een{\]}
\def\bearn{\begin{eqnarray*}}
\def\eearn{\end{eqnarray*}}
\def\bear{\begin{eqnarray}}
\def\eear{\end{eqnarray}}
\def\barr{\begin{array}}
\def\earr{\end{array}}
\def\bmat{\left(\begin{array}}
\def\emat{\end{array}\right)}

\newcommand{\limit}[2]{\lim_{#1 \rightarrow #2}}
\newcommand{\eqdef}{\stackrel{\mathrm{def}}{=}} % by definition

\newcommand{\dert}[1]{\frac{d#1}{dt}} % derivative wr t
\newcommand{\mand} {\mathrm{\;  and \; }}
\newcommand{\mor} {\mathrm{\;  or \; }}

\newcommand{\mfa} {\mathrm{\;  for \;  all \; }}

\def\Reals{\mathbb{R}}

\def\Nats{\mathbb{N}}
\def\E{\mathbb{E}}
\def\P{\mathbb{P}}

\newcommand{\bracket}[1]{ \left\{\begin{array}{l} #1 \end{array} \right.}

\newsavebox{\coloredbox}

\newsavebox{\traitbox}

\title{ 
\begin{minipage}{9cm}
\scriptsize{
\textsf{Post-publication version of\\~\\
http://dx.doi.org/10.3934/nhm.2013.8.529}, \\ 
NETWORKS AND HETEROGENEOUS MEDIA\\ Volume 8, Number 2, June 2014, pp. 529-540.
~\\}
\end{minipage}
\\~\vspace{2cm}\\
The Stationary Behaviour of Fluid Limits of Reversible Processes is Concentrated on Stationary Points}
\author{Jean-Yves Le Boudec%
\footnote{
  EPFL IC-LCA2 - Lausanne, Switzerland.\\
  \textsf{\textbf{Erratum: } In the publication version, the definition of $f_i^N$, just after \eqref{eq-sp}, should be $f_i^N(n) \eqdef \frac{\theta_i^n}{n!}\frac{1}{\prod_{m=1}^n \mu^N(m)}$ instead of $f_i^N(n) \eqdef \frac{\theta_i^n}{n!}\prod_{m=1}^n \mu^N(m)$. This is fixed in this version and is the only difference with the publication version.
}  }
  }
\begin{document}

%\maketitle
\begin{frontmatter}

\begin{abstract}
  Assume that a stochastic process can be approximated, when some scale parameter gets large, by a fluid limit (also called ``mean field limit", or ``hydrodynamic limit"). A common practice, often called the ``fixed point approximation" consists in approximating the stationary behaviour of the  stochastic process by the stationary points of the fluid limit. It is known that this may be incorrect in general, as the stationary behaviour of the fluid limit may not be described by its stationary points. We show however that, if the stochastic process is reversible, the fixed point approximation is indeed valid. More precisely, we assume that the stochastic process converges to the fluid limit in distribution (hence in probability) at every fixed point in
time. This assumption is very weak and holds for a large family of processes,
among which many mean field and other interaction models. We show that the reversibility of the stochastic process implies that any limit point of its stationary distribution is concentrated on stationary points of the fluid limit. If the fluid limit has a unique stationary point, it is an approximation of the stationary distribution of the stochastic process. 
\end{abstract}

\end{frontmatter}
%\section*{AMS Subject Classification}
%37 Dynamical Systems and Ergodic Theory
%
%60 Probability Theory and Stochastic Processes
\newpage
\section{Introduction}
%
%[tbd: It seems that we do not need to assume that the semi-flow
%is continuous for Theorem 1 to hold].

This paper is motivated by the use of fluid limits in models of interacting objects or particles, in contexts such as communication and computer system modelling \cite{prout07}, biology \cite{borghans1998tcellvaccination} or game theory \cite{benaim03}.
%
%common use of the fixed point approximation in models of interacting objects, as found in such contexts as communication and computer system modelling \cite{kelly1991loss,bianchi,altman07,LCA-ARTICLE-2008-034}, biology \cite{borghans1998tcellvaccination} or game theory \cite{}.
Typically, one has a  stochastic process $Y^N$, indexed by a size parameter $N$; under fairly general assumptions, one can show that the stochastic process $Y^N$ converges to a deterministic fluid limit $\varphi$ \cite{kurtz-81}. We are interested in the stationary distribution of $Y^N$, assumed to exist and be unique, but which may be too complicated to be computed explicitly. The ``fixed point assumption" is then sometimes invoked \cite{altman07,bianchi,LCA-CONF-2006-015,kelly1991loss}: %\cite{LCA-ARTICLE-2008-034}
it consists in approximating the stationary distribution of $Y^N$ by a stationary point of the deterministic fluid limit $\varphi$. In the frequent case where the fluid limit $\varphi$ is described by an Ordinary Differential Equation (ODE), say of the form $\dot{y}=F(y)$, the stationary points are obtained by solving $F(y)=0$.
If $Y^N$ is an empirical measure, convergence to a deterministic limit implies propagation of chaos, i.e. the states of different objects are asymptotically independent, and the distribution of any particular object at any time is obtained from the fluid limit. Under the fixed point assumption, the stationary distribution of one object is approximated by a stationary point of the fluid limit.

A critique of the fixed point approximation method is formulated in \cite{benaim2008class}, which observes that one may only say, in general, that the stationary distribution of $Y^N$ converges to a stationary distribution of the fluid limit. For a deterministic fluid limit, a stationary distribution is supported by the Birkhoff center of the fluid limit, which may be larger than the set of stationary points. An example is given  in \cite{benaim2008class} where the fluid limit has a unique stationary point, but the stationary distribution of $Y^N$ does not converge to the Dirac mass at this stationary point; in contrast, it converges to a distribution supported by a limit cycle of the ODE. If the fluid limit has a unique limit point, say $y^*$, to which all trajectories converge, then this unique limit point is also the unique stationary point and the stationary distribution of $Y^N$ does converge to the Dirac mass at $y^*$ (i.e. the fixed point approximation is then valid). However, as illustrated in \cite{benaim2008class}, this assumption may be difficult to verify, as it often does not hold, and when it does, it may be difficult to establish. For example, in \cite{EPFL-CONF-149779} it is shown that the fixed point assumption does not hold for some parameter settings of a wireless system analyzed in \cite{bianchi}, due to limit cycles in the fluid limit.

In this paper we show that there is a class of systems for which such complications may not arise, namely the class of reversible stochastic processes. Reversibility is classically defined as a property of time reversibility in stationary regime \cite{kelly1979reversibility}. There is a large class of processes that are known to be reversible, for example product-form queuing networks with reversible routing, or stochastic processes in \cite{kelly1991loss}, which describes the occupancy of inter-city telecommunication links; in \sref{sec-ex} we give an example motivated by crowd dynamics.
 In such cases, we show that the fluid limit must have stationary points, and any limit point of the stationary distribution of $Y^N$ must be supported by the set of stationary points. Thus, for reversible processes that have a fluid limit, the fixed point approximation is justified. 

\section{Assumptions and Notation}
\label{assu}

\subsection{A Collection of Reversible Random Processes} Let $E$ be a Polish space and let $d$ be a measure that metrizes $E$. Let $\calP(E)$ be the set
of probability measures on $E$, endowed with the topology of
weak convergence. Let $\calC_b(E)$ be the set of bounded
continuous
functions from $E$ to $\Reals$, and similarly $\calC_b(E\times E)$ is the set of bounded
continuous
functions from $E\times E$ to $\Reals$.

We are given a collection of probability spaces
$(\Omega^N,\calF^N, \P^N)$ indexed by $N =1,2,3,...$ and for
every $N$ we have a process $Y^N$ defined on
$(\Omega^N,\calF^N, \P^N)$. Time is
continuous. Let $D_E[0,\infty)$ be the set of
c\'adl\'ag functions $[0,\infty) \to E$; $Y^N$ is then a
stochastic process with sample paths in $D_{E}[0,\infty)$.

We denote by $Y^N(t)$ the random value of $Y^N$ at time
$t\geq 0$. Let $E^N\subset E$ be the support of $Y^N(0)$, so
that $\P^N(Y^N(0)\in E^N)=1$.

We assume that, for every $N$, the process $Y^N$ is Feller, in
the sense that for every $t\geq 0$ and $h \in \calC_b(E)$,
$\E^N\left. \lb h(Y^N(t))\right| Y^N(0)=y_0\rb$ is a
continuous function of $y_0 \in E$. Examples of such processes are continuous time Markov chains as in \cite{kurtz1970solutions}, or linear interpolations of discrete time Markov chains as in \cite{benaim-07}, or the projections of a Markov process as in \cite{graham1993propagation}.
%Note that apart from the first example, these are not Markov.

\begin{defi} A probability $\Pi^N\in \calP(E)$ is \emph{invariant} for $Y^N$ if
$\Pi^N(E^N)=1$ and for every $h \in \calC_b(E)$ and every $t
\geq 0$:
  \ben
 \int_E \E^N\left. \lb h\lp Y^N(t)\rp\right|Y^N(0)=y\rb \Pi^N(dy) = \int_E h(y)
 \Pi^N(dy)
  \een\label{def1}
  \end{defi}

We are interested in reversible processes, i.e. processes that keep the same stationary law under time reversal. A weak form of such a property is defined as follows

\begin{defi} Assume $\Pi^N$ is a probability on $E$ such that $\Pi^N(E^N)=1$, for some $N$. We say that $Y^N$ is \emph{reversible} under $\Pi^N$ if for every time $t \geq 0$ and any $h \in \calC_b(E\times E)$:
 \ben
 \int_E \E^N\left. \lb h\lp y, Y^N(t)\rp\right|Y^N(0)=y\rb \Pi^N(dy)
 =
 \int_E \E^N\left. \lb h\lp Y^N(t),y \rp\right|Y^N(0)=y\rb \Pi^N(dy)
 \een
\label{defr}
\end{defi}
Note that, necessarily, $\Pi^N$ is an invariant probability for $Y^N$.
If $Y^N$ is an ergodic Markov process with enumerable state space, then Definition~\ref{defr} coincides with the classical definition of reversibility by Kelly in \cite{kelly1979reversibility}. Similarly, if $Y^N$ is a projection of a reversible Markov process $X^N$, as in \cite{crametz1991limit}, then $Y^N$ is reversible under the projection of the stationary probability of $X^N$; note that in such a case, $Y^N$ is not Markov.

\subsection{A Limiting, Continuous Semi-Flow}

Further, let $\varphi$ be a deterministic process, i.e. a
measurable mapping
 \ben \barr{rccl}
 \varphi:& [0,\infty) \times E &\to &E
 \\
&t,y_0 & \mapsto & \varphi_t(y_0)
 \earr\een

We assume that $\varphi_t$ is a semi-flow, i.e.
\begin{enumerate}
\item $\varphi_0(y) = y$,
     \item $\varphi_{s+t}=\varphi_s\circ
\varphi_t$ for all $s\geq 0$ and $t\geq 0$,
   \end{enumerate}
and we say that $\varphi$ is ``space continuous" if for every $t \geq 0$, $\varphi_{t}(y)$ is continuous in $y$.

\begin{defi}
We say that $y \in E$ is a \emph{stationary point} of $\varphi$ if $\varphi_t(y)=y$ for all $t \geq 0$
\end{defi}
In cases where $E$ is a subset of $\Reals^d$ for some integer $d$, the semi-flow $\varphi$ may be an autonomous ODE, of the form $\dot{y}=F(y)$; here the stationary points are the solutions of $F(y)=0$.

%Similar to Definition~\ref{defr}, we have:
 \begin{defi} We say that the semi-flow $\varphi$ is \emph{reversible} under the probability $\Pi \in \calP(E)$ if for every time $t \geq 0$ and any $h \in \calC_b(E\times E)$:
  \be
\int_E h(y, \varphi_t(y)) \Pi(dy)= \int_E h(\varphi_t(y),y) \Pi(dy)
\ee
 \label{defrs}
 \end{defi}  As we show in the next section, reversible semi-flows must concentrate on stationary points.

%be a continuous semi-flow defined on $E$, i.e. for every $t$,
%$\varphi_t(y_0)$ is a continuous function of $y_0$.

%
%Examples: ODEs with unique limit point; with limit circle.

\subsection{Convergence Hypothesis} We assume that, for every fixed $t$ the processes
$Y^N$ converge in distribution to some space continuous deterministic process
$\varphi$ as $N \to \infty$ for every collection of converging
initial conditions. More precisely:

\begin{assh}\label{ass-1}For every $y_0$ in $E$, every sequence
$(y^N_0)_{N=1,2,...}$ such that $y^N_0  \in E^N$ and
$\limit{N}{\infty}y^N_0=y_0$, and every $t\geq0$, the
conditional law of $Y^N(t)$ given $Y^N(0)=y^N_0$ converges in distribution to the Dirac mass at $\varphi_t(y_0)$. That is
\be\limit{N}{\infty} \E^N\left. \lb h(Y^N(t))\right| Y^N(0)=y_0^N\rb = h \circ \varphi_t(y_0)
\label{eq-h1}
\ee for all
$h \in \calC_b(E)$ and any fixed $t\geq 0$. In the above, $\varphi$ is a space continuous semi-flow.
\end{assh}

Hypothesis~\ref{ass-1} is commonly true in the context of fluid or mean field limits. 
The stronger convergence
results results in
\cite{le2007generic,kurtz1970solutions,sandholm2006population,prout07} imply that Hypothesis~\ref{ass-1} is satisfied ; we give a detailed example illustrating this in \sref{sec-ex}.
Similarly, \cite{benaim2008class} gives very general conditions (called H1 to H5) that ensure convergence of a stochastic process to its mean field limit; under these conditions,  Hypothesis~\ref{ass-1} is automatically satisfied  (the deterministic process $\varphi$ is then an ODE).
Note that the results in these references are stronger than what we require in Hypothesis~\ref{ass-1}; for example in \cite{kurtz1970solutions} there is
almost sure, uniform convergence for all $t \in [0,T]$, for any
$T\geq 0$; in \cite{graham1993propagation} the convergence is on the set of trajectories.

Under Hypothesis~\ref{ass-1}, $\varphi$ is called the \emph{hydrodynamic limit}
%\cite{aldous1999deterministic},
or simply \emph{fluid limit} of $Y^N$.

%
%AJOUTER LOSS NETWORK ET GRAHAM

\section{Reversible Semi-Flows Concentrate on Stationary Points}
\begin{thm} Let $\varphi$ be a space continuous semi-flow, reversible under $\Pi$. Let $S$ be the set of stationary points of $\varphi$. Then $\Pi$ is concentrated on $S$, i.e. $\Pi(S)=1$.
\label{theo1}
\end{thm}
\begin{mypr}

  \textbf{Step 1. }
    Denote with $S^c$ the complement of the set of stationary points. Take some fixed but arbitrary $y_0\in S^c$. By definition of $S$, there exists some $\tau>0$ such that
 \be
 \varphi_{2 \tau}(y_0) \neq y_0
 \ee
Define $ \varphi_{\tau}(y_0)  =  y_1$, $\varphi_{\tau}(y_1)  =  y_2$,   so that $y_2  \neq  y_0$.

For $y \in E$ and $\eps >0$ we denote with $B(y,\eps)$ the open ball $=\lc x \in E, d(x,y) < \eps \rc$.
Let $\eps= d(y_0,y_2)>0$ and
let $B_2=B(y_2, \eps/2)$.
Since the semi-flow is space continuous, there is some $\alpha_1 >0$ such that $B_1=B(y_1, \alpha_1)$ and $\varphi_{\tau}\lp B_1 \rp  \subset  B_2$. Also let $B_1'=B(y_1, \alpha_1/2)$. By the same argument, there exists some $\alpha_0>0$ such that $\alpha_0<\eps/2$, $B_0=B(y_0, \alpha_0)$ and $\varphi_{\tau}\lp B_0 \rp  \subset  B_1'$.
We have thus:
  \bearn
   \varphi_{\tau}\lp B_0 \rp & \subset & B'_1 \subset B_1
   \\
   \varphi_{\tau}\lp B_1 \rp & \subset & B_2
   \\
   B_0 \cap B_2 & = & \O
  \eearn
Let $\xi$ be some continuous function $[0,+\infty) \to [0,1]$ such that $\xi(u)=1$ whenever $0\leq u \leq 1/2$ and $\xi(u)=0$ whenever $u \geq 1$ (for example take a linear interpolation).
Now take
 \be
 h(y,z)\eqdef \xi\lp \frac{ d(y_0,y)}{\alpha_0}\rp  \xi\lp \frac{ d(y_1,y)}{\alpha_1}\rp
 \ee
 so that $h \in \calC_b(E\times E)$ and
  \bearn
   h(y,z) = 0 && \mbox{ whenever } y \nin B_0 \mor z \nin B_1
   \\
   h(y,z) = 1 && \mbox{ whenever } d(y_0,y)<\alpha_0/2 \mand  z \in B'_1
  \eearn
  It follows that
   $h(\varphi_{\tau}(z),z)= 0$ for every $z \in E$ and
     \be
     \int_E  h(y,\varphi_{\tau}(y)) \Pi(dy) \geq \Pi(B(y_0,\alpha_0/2))
     \ee
    Apply Definition~\ref{defrs}, it comes $\Pi\lp B(y_0,\alpha_0/2) \rp = 0$;
    thus, for any non stationary point $y_0$ there is some $\alpha >0$ such that
  \be
  \Pi\lp B(y_0,\alpha) \rp = 0 \label{eq-rev-s2}
  \ee

%For $y \in E$ and $\eps >0$ we denote with $B(y,\eps)$ the open ball $=\lc x \in E, d(x,y) < \eps \rc$.
%Let $\eps= d(y_0,y_2)>0$ and
%let $B_2=B(y_2, \eps/2)$. Also let $B_1= \varphi_{\tau}^{-1}(B_2)$. Since the semi-flow is continuous in space, $B_1$ is a neighbourhood of $y_1$ and $\varphi_{\tau}^{-1}(B_1)$ is a neighbourhood of $y_0$. Thus there exists some $\alpha < \eps/2$ such that the ball $B_0=B(y_0,\alpha)$ is included in $\varphi_{\tau}^{-1}(B_1)$. We have thus:
%  \bearn
%   \varphi_{\tau}\lp B_0 \rp & \subset & B_1
%   \\
%   \varphi_{\tau}\lp B_1 \rp & \subset & B_2
%   \\
%   B_0 \cap B_2 & = & \O
%  \eearn
%
%Now take
% \be
% h(y,z)\eqdef \ind{y\in B_0}\ind{z \in B_1}
% \ee It follows that $h(y, \varphi_{\tau}(y))= \ind{y \in B_0}$ and $h(\varphi_{\tau}(y),y)= 0$. Apply Definition~\ref{defrs}, it comes:
%  \be
%  \Pi\lp B(y_0,\alpha) \rp = 0 \label{eq-rev-s2}
%  \ee
% Thus, for any non stationary point $y_0$ there is some $\alpha >0$ such that \eref{eq-rev-s2} holds.

\textbf{Step 2. }
The space is polish thus also separable, i.e. has a dense enumerable set, say $Q$.

 For every $y \in S^c$ let $\alpha$ be as in \eref{eq-rev-s2} and pick some $q(y) \in Q$ and $n(y) \in \Nats$ s.t. $d(y,q(y))< \frac{1}{n(y)}<\alpha$. Thus $y \in B(q(y)$, $\frac{1}{n(y)})$ and $\Pi\lp B(q(y),\frac{1}{n(y)} )\rp =0$.

 Let $F=\bigcup_{y \in S^c}(q(y),n(y))$. $F\subset Q \times \Nats$ thus $F$ is enumerable and
   \bearn
   S^c &\subset & \bigcup_{(q, n) \in F} B\lp q, \frac{1}{n}\rp   \eearn
   Thus
   \be
   0\leq \Pi(S^c) \leq \sum_{(q, n) \in F} \Pi\lp B\lp q, \frac{1}{n}\rp \rp =0
   \ee
\end{mypr}

Note that it follows that a semi-flow that does not have any stationary point cannot be reversible under any probability.

\section{Stationary Behaviour of Fluid Limits of Reversible Processes}

\begin{thm}Assume for every $N$ the process $Y^N$ is reversible under some probability $\Pi^N$. Assume the convergence Hypothesis~\ref{ass-1} holds and that $\Pi \in \calP(E)$ is a limit point (for weak convergence) of the
sequence $\Pi^N$. Then the fluid limit is reversible under $\Pi$. In particular, it follows from \thref{theo1} that $\Pi$ is concentrated on the set of stationary points $S$ of the fluid limit $\varphi$.
\label{theo2}
\end{thm}
\begin{mypr}
All we need to show is that $\Pi$ verifies Definition~\ref{defrs}.  %Using Skorohod's representation theorem for Polish
%spaces \cite[Thm 1.8]{ethier-kurtz-05}, dominated convergence as in \cite{lebbenaimstat2010}, and Definition~\ref{defr}, we obtain:
  Let $N_k$ be a subsequence such that
$\limit{k}{\infty}\Pi^{N_k}=\Pi$ in the weak topology on
$\calP(E)$. By Skorohod's representation theorem for Polish
spaces \cite[Thm 1.8]{ethier-kurtz-05}, there exists a common
probability space $(\Omega, \calF,\P)$ on which some random
variables $X^k$ for $k \in \Nats$ and $X$ are defined such that
 \ben
 \bracket{
 \mbox{ law of } X^{k} = \Pi^{N_k}\\
 \mbox{ law of } X = \Pi\\
 X^{k}\to X \; \P-\mbox{a.s.}
  }
 \een

Fix some $t \geq 0$ and $h \in \calC_b(E\times E)$, and define, for $k
\in \Nats$ and $y \in E$
  \bearn
  a^k(y) &\eqdef& \espc{h\lp y, Y^{N_k}(t)\rp}{Y^{N_k}(0)=y}
  \\
  b^k(y) &\eqdef& \espc{h\lp Y^{N_k}(t),y\rp}{Y^{N_k}(0)=y}
    \eearn

Since $Y^N$ is reversible under $\Pi^{N_k}$:
  \be
  \int_E a^k(y) \Pi^{N_k}(dy) = \int_E b^k(y) \Pi^{N_k}(dy)\label{eq-877}
  \ee
 %which we can also write as
% \be
% \esp{a^k(X^k)} = \esp{h(X^k)}
% \ee

Hypothesis~\ref{ass-1} implies that
$\limit{k}{\infty}a^k(x^k)=h(x,\varphi_t(x))$ for every sequence
$x^k$ such that $x^k\in E^{N_k}$ and $\limit{k}{\infty}x^k=x
\in E$. Now $X^k \in E^{N_k}$ $\P-$ almost surely, since the law
of $X^k$ is $\Pi^{N_k}$ and
$Y^{N_k}$ is reversible under $\Pi^{N_k}$ . Further, $X^k\to X$ $\P-$ almost surely; thus

 \be
  \limit{k}{\infty} a^k(X^k) = h(X,\varphi_t(X)) \;\;\;
  \P-\mbox{ almost surely}
  \ee
 Now $a^k(X^k)\leq \norm{h}_{\infty}$ and, thus, by dominated
 convergence:
 \be
 \limit{k}{\infty}\esp{a^k(X^k)}=\esp{h(X,\varphi_t(X))}
 \ee
and similarly for $b^k$. Thus
    \be
\int_E h(y, \varphi_t(y)) \Pi(dy)= \int_E h(\varphi_t(y),y) \Pi(dy) \label{eq-rev-s33}
\ee
\end{mypr}

 In particular, if the semi-flow has a
 unique stationary point, we have:

 \begin{coro}\label{coro-uni}
 Assume the processes $Y^N$ are reversible under some probabilities $\Pi^N$. Assume Hypothesis~\ref{ass-1} holds and:
  \begin{enumerate}
    \item the sequence $(\Pi^N)_{N=1,2,...}$ is tight;
    \item the semi-flow $\varphi$ has a unique stationary point $y^*$.
  \end{enumerate}
  It follows that the sequence $\Pi^N$ converges weakly to the Dirac mass at
  $y^*$.
 \end{coro}
We leave the proof  of the corollary to the reader (it follows in a classical way from compactness arguments; the tightness condition implies that the set $\lc \Pi^N, N=1,2,3... \rc_N$ is relatively compact in $\calP(E)$).

Recall that tightness means that for every $\eps>0$ there is
 some compact set $K \subset E$ such that $\Pi^N(K)\geq 1-\eps$
 for all $N$. If $E$ is compact then $(\Pi^N)_{N=1,2,...}$ is
 necessarily tight, therefore condition 1 in the corollary is automatically satisfied. For mean field limits where $E$ is the simplex in finite dimension, the corollary says that, if the pre-limit process is reversible, then the existence of a unique stationary point implies that the Dirac mass at this stationary point is the limit of the stationary probability of the pre-limit process.

 Compare \coref{coro-uni} to known results for the non reversible case \cite{benaim98}: there we need that the fluid limit $\varphi$ has a unique limit point to which all trajectories converge. In contrast, here, we need a much weaker assumption, namely, the existence and uniqueness of a stationary point. It is possible for a semi-flow to have a unique stationary point, without this stationary point being a limit of all trajectories (for example because it is unstable, or because there are stable limit cycles as in \cite{benaim2008class}). In the reversible case, we do not need to show stability of the unique stationary point $y^*$. 
\section{Example: Crowd Dynamics}
\label{sec-ex}
In this section we give an example to illustrate the application of \thref{theo2} -- a detailed study of this example beyond the application of \thref{theo2} is outside the scope of this paper. We consider the crowd dynamics model of \cite{RoG03}. The model captures the emergence of crowds in a city. A city is modelled as a fully connected bidirectional graph with $I$ vertices, every vertex representing a square, where bars are located. There is a fixed total population $N$. People spend some time in a square and once in a while decide to leave a square and move to some other square. The original model is in discrete time, and at every time slot, the probability that a tagged person present in square $i$ leaves square $i$ is assume to be equal to
$(1-c)^{N_i(t)-1}$ where $N_i(t)$ is the population of square $i$ at time $t$. In this equation, $c$ is the \emph{chat} probability, and this model thus assumes that a person leaves a square when it has no one to chat with. The model also assumes that departure events are independent. When a person leaves a square $i$, it move to some other square $j$ according to Markov routing, with probability $Q_{i,j}$ given by
\be
 Q_{i,j} = \frac{1}{d(i)}
\label{eq-Q}
\ee where $d(i)$ is the degree of node $i$, i.e. the person picks a neighboring square $j$ uniformly at random among all neighboring squares.

In \cite{RoG03}, the authors study by simulation the emergence of concentration in one square. They also show that for regular graphs (i.e. when all vertices have same degree) there is a critical value
%$c^*=I/N$
$c^*$
such that for $c>c^*$ concentrations occur, whereas for $c<c^*$ the stationary distribution of people is uniform. The analysis is based on the study of stationary points for the empirical distribution. Note that, as mentioned in the introduction, the analysis with stationary points may, in general, miss the main part of the stationary distribution, and it is quite possible that the stationary distribution is not concentrated on stationary points (for example if there is a limit cycle \cite{benaim2008class}). A fluid flow approximation is proposed in \cite{Ma+11}, and similar results are found.

To understand whether the stationary point analysis is justified, we study the large $N$ asymptotics for an appropriately rescaled version. To avoid unnecessary complications, we replace the original model of \cite{RoG03}, which is in discrete time, by its continuous time counterpart. The probability that a tagged person present in square $i$ leaves square $i$ in a time slot is replaced by the rate of service given by
\be
\mu^N(N_i) \eqdef  (1-c)^{N_i-1}
\ee
i.e. the probability that a tagged person present in square $i$ leaves the square in the next $dt$ seconds is $\mu^N(N_i)dt + o(dt)$.

%&where $\Delta_N$ is the duration in continuous time of one discrete time slot.
The corresponding continuous process $X^N(t)=(N_1(t), ...,N_I(t))$ is a Markov process on an enumerable state space. More precisely, it is a queuing network of infinite server stations, with state dependent service rate and with Markov routing. It follows from classical results on quasi-reversibility that it has product-form (see for example \cite[Chapter 8]{leboudec2010performance}), i.e. it is ergodic (since the graph of squares is fully connected and the population is finite) and its stationary probability is given, for every $(n_1, ...,n_I)\in \Nats^I$ such that $n_1+...+n_I=N$ by
\bear
 P^N(n_1,...,n_I)&=&\eta^N \prod_{i=1}^I f^N_i(n_i)  
 \label{eq-sp}\eear In this formula, $\eta^N$ is a normalizing constant, $f_i^N(n) \eqdef \frac{\theta_i^n}{n!}\prod_{m=1}^n \mu^N(m)$, and the vector $\theta=(\theta_1,...,\theta_I)$ is the stationary distribution of the Markov routing matrix $Q$, i.e. a normalized solution of the equation $\theta Q=\theta$. Note that it follows from \eref{eq-Q}
that
\be
\theta_i=\frac{d(i)}{\sum_{j=1}^I d(j)}
\label{eq-theta}
\ee
%The stationary probability of $X^N$ has a closed form, however, in practice, the computation can only be performed numerically, and it is difficult to say anything general from it.

To apply \thref{theo2}, we need to show that $X^N$ is reversible.
A product-form queuing network is, in general, not reversible. However, it is so if the Markov routing chain is reversible \cite{le1987interinput}, which is the case here.
\begin{thm}For every $N$, the process $X^N$ is reversible.
\end{thm}
\begin{mypr} Take $\theta$ given by \eref{eq-theta}. Then $\theta_i Q_{i,j}=\theta_j Q_{j,i}$ for any pair $(i,j)$, thus the Markov chain with transition matrix $Q$ given by \eref{eq-Q} is reversible. By \cite{le1987interinput}, it follows that the product-form queuing network $X^N$ is reversible.
\end{mypr}

In \cite{bortolussi2012revisiting}, it is suggested to scale the chat probability as
\be c=\frac{s}{N}\label{eq-c}\ee in order to account for the fact that, for large populations, meetings tend to be limited by space or size of the friend's group. We use this scaling law and %take as time slot duration $\Delta_N=\frac{1}{N}$. We then
consider the re-scaled process $Y^N$ of occupancy measures, i.e. given by
 \be
 Y^N(t)=\lp\frac{N_1(t)}{N},...\frac{N_I(t)}{N}\rp = \frac{1}{N}X^N(t)
 \ee Obviously, for every $N$ the process $Y^N$ is Markov and is reversible. Further, it converges to an ODE, as we see next.

 To establish the convergence of $Y^N$, we compute its drift
 \be
 A^N(y)\eqdef \limit{dt}{0}\frac{\espc{Y^N(t+dt)-y}{Y^N(t)=y}}{dt}
 \ee
defined for every possible value $y$ of $Y^N(t)$. When the occupancy measure is $y$, there are $N_i=Ny_i$ people in square $i$, and the rate of departure from square $i$ is $N_i \mu(N_i)=Ny_i \mu(N y_i)$; the delta to the occupancy measure due to one person moving from square $i$ to square $j$ is $\frac{-e_i+e_j}{N}$, where $e_i$ is the row vector with a $1$ in position $i$ and $0$ elsewhere. Therefore
\ben
V^N(y)=
\sum_{i,j} N y_i \mu(Ny_i)Q_{i,j} \frac{-e_i+e_j}{N}
%\sum_{i,j} N y_i (1-c)^{N y_i-1}Q_{i,j} \lp-\ve_i+\ve_j \rp
\een
Taking into account \eref{eq-c}, it comes
\be
V^N(y)= \sum_{i,j} y_i \lp 1- \frac{s}{N}\rp^{Ny_i-1}Q_{i,j} \lp -e_i+e_j\rp
\ee
Let $\Delta_I$ denote the simplex, i.e.
\ben
\Delta_I=\lc y\in\Reals^I, y_i \geq 0 \mfa i \mand \sum_{i=1}^I y_i=1 \rc
\een
 When $N\to \infty$, $V^N(y)$ converges for every $y\in \Delta_I$ to
\be
V(y) \eqdef \sum_{i,j} y_i e^{-sy_i}Q_{i,j} \lp -e_i+e_j\rp
\ee
This suggests that the fluid limit of $Y^N$, if it exists, would be the deterministic process $y(t)$, with sample paths in $\Reals^d$, obtained as solution of the ODE $\dert{y}=V(y)$. As we show next, this is indeed the case and follows from ``Kurtz's theorem" \cite[Theorem 9.2.1]{sandholm2006population}. Before that, we rewrite the ODE more explicitly as
\be
\dert{y_i}=-y_i e^{-sy_i} + \sum_j y_j e^{-sy_j}Q_{j,i}
\label{eq-ode}
\ee
\begin{thm}Assume that the initial condition $y_0^N$ of the process $Y^N$ is deterministic and converges to some $y_0\in\Delta_I$. Let $\varphi$ be the semi-flow defined by the ODE~(\ref{eq-ode}), i.e. $\varphi_t(y_0)$ is the solution of the ODE~(\ref{eq-ode}) with initial condition $y(0)=y_0$ (this solution exists and is unique by the Cauchy Lipschitz theorem). Then for each $T>0$ and $\eps>0$:
\be
\limit{N}{\infty}\P\lp\sup_{0\leq t \leq T}\norm{Y^N(t)-\varphi_t(y_0)}>\eps\rp=0
\label{eq-conv}
\ee
(the notation $\norm{}$ stands for any norm on $\Reals^I$).

It follows that Hypothesis 1 is verified.
\label{theo4}
\end{thm}
\begin{mypr}
We apply Theorem 9.2.1 in \cite{sandholm2006population}. We need to find a sequence of numbers $\delta_N\to 0$ such that the following three conditions hold:
\bear
&&\limit{N}{\infty}\sup_{y\in \Delta_I^N}\norm{V^N(y)-V(y)}=0\label{eq-cond1}\\
&&\sup_N \sup_{y\in \Delta^N_I}A^N(y)<\infty\label{eq-cond2}\\
&&\limit{N}{\infty} \sup_{y\in \Delta^N_I}A^N_{\delta_N}(y)=0\label{eq-cond3}
\eear In the above, $\Delta^N_I$ is the set of feasible states of $Y^N$, i.e. the set of  $y\in \Delta_I$ such that $Ny$ is integer, $A^N(y)$ is the expected norm of jump per time unit, and $A^N_{\delta_N}(y)$ is the absolute expected norm of jump per time unit due to jumps travelling further than $\delta_N$.

We now show \eref{eq-cond1}. First consider the case $y\in \Delta^N_I$ such that $y_i>0$ (thus we have $1/N\leq y_i\leq 1$). We apply the inequality
$\abs{e^{-x} -e^{-x'}}\leq \abs{x-x'}$, valid for $x\geq0$ and $x'\geq0$, to $x=-(Ny_i-1)\log\lp 1-\frac{s}{N}\rp$, $x'=sy_i$ and obtain:
\ben
 \abs{\lp 1-\frac{s}{N}\rp^{Ny_i -1}-e^{-s y_i}}\leq \abs{(Ny_i-1)\log\lp 1-\frac{s}{N}\rp +sy_i}
\een The right handside is convex in $y_i$ thus its maximum for $y_i\in [1/N, 1$] is obtained at one end of the interval. Thus
\ben
 \abs{\lp 1-\frac{s}{N}\rp^{Ny_i -1}-e^{-s y_i}}\leq a_N(s)
 \eqdef \max\lb {s\over N}, \abs{(N-1) \log\lp 1-\frac{s}{N}\rp +s}\rb
\een
Second, multiply by $y_i$ and note that $y_i\leq 1$, it follows that, whenever $y\in \Delta^N_I$ and $y_i>0$:
\ben
 \abs{y_i \lp 1-\frac{s}{N}\rp^{Ny_i -1}-y_i e^{-s y_i}}\leq a_N(s)
\een and this is also obviously true if $y_i=0$. It follows that
\ben
\sup_{y\in\Delta^N_I}\norm{V^N(y)-V(y)}\leq a_N(s) \sum_{i,j}\norm{-e_i+e_j}Q_{i,j}
\een from where \eref{eq-cond1} follows since $\limit{N}{\infty}a_N(s)=0$.

We now show \eref{eq-cond2}. We take the sup norm on $\Reals^d$ so that $\norm{-e_i+e_j}=1$ for $i\neq j$; thus
\bearn
A^N(y)&\eqdef& \sum_{i,j} y_i \lp 1- \frac{s}{N}\rp^{Ny_i-1}Q_{i,j} \norm{ -e_i+e_j}
\\
&=&
\sum_{i,j} y_i \lp 1- \frac{s}{N}\rp^{Ny_i-1}Q_{i,j}\leq \sum_{i,j}Q_{i,j} = I
\eearn which shows \eref{eq-cond2}.

We now show \eref{eq-cond3}. We take $\delta_N=\frac{1}{N}$. The jumps of $Y^N$ are of the form $\frac{-e_i+e_j}{N}$ and thus $\norm{\frac{-e_i+e_j}{N}}\leq\delta_N$. Thus
\ben
A^N(y)_{\delta_N}\eqdef \sum_{i,j} y_i \lp 1- \frac{s}{N}\rp^{Ny_i-1}Q_{i,j} \norm{ -e_i+e_j}\ind{\norm{\frac{-e_i+e_j}{N}}>\delta_N}=0
\een which trivially shows \eref{eq-cond3}. By Theorem 9.2.1 in \cite{sandholm2006population}, this establishes \eref{eq-conv}.

It follows (as a much weaker convergence) that for any fixed $T$, $Y^N(T)$ converges in probability to the deterministic $\varphi_T(y_0)$. Thus there is also convergence in law, since convergence in probability to a deterministic variable implies convergence in distribution, i.e. \eref{eq-h1} in Hypothesis 1 is verified. It remains to see that $\varphi_t$ is space continuous: this follows from the fact that the right-handside of the ODE is Lipschitz continuous and from the Cauchy Lipschitz theorem.
\end{mypr}

It follows from \thref{theo2} that any limit point of the stationary probability \eref{eq-sp} is concentrated on the stationary points of the ODE (\ref{eq-ode}). This justifies a posteriori the method in \cite{RoG03}, which looked only at stationary points.

For the case of a regular graph (this is the case studied analytically in \cite{RoG03}), the stationary points can be obtained explicitly (Theorem 6.1 in  \cite{bortolussi2012revisiting}). In particular, there is a critical value $s^*$ below which there is only one stationary point, equal to the uniform distribution $y^*=(\frac{1}{I},...,\frac{1}{I})$ and above which there are other stationary points. The critical value is given in \cite{bortolussi2012revisiting}) and is equal to
\be
s^*= \min_{K=1,...,I-1}\min_{\alpha>1}\lp (I-K)\alpha + K \phi(\alpha)\rp
\ee  with $\phi(x)\eqdef -W_0(-x e^{-x})$, $W_0$ being the Lambert-W function of index $0$. For example for $I=3$, the critical value is $s^*\approx 2.7456$.

We can apply \coref{coro-uni}: since the state space $E=\Delta_I$ is compact, it follows that for $s<s^*$, the stationary distribution given by \eref{eq-sp}, re-scaled by $1/N$, converges as $N\to \infty$ to the uniform distribution. This illustrates the interest of the reversibility results in this paper; we do not need to show that all trajectories converge to the single stationary point -- its uniqueness and the reversibility argument are sufficient. 

\bibliographystyle{plain}
\bibliography{leb}

\end{document}